\documentclass[11pt]{elsarticle}
\usepackage[all,pdf]{xy}

\usepackage{amsmath,amssymb,amsthm}
\usepackage{hyperref}
\usepackage{mathrsfs}
\usepackage{graphicx}
\usepackage[all,pdf]{xy}
\usepackage{tikz}
\usepackage{titlesec}

\theoremstyle{plain}
\newtheorem{theorem}{Theorem}[section]
\newtheorem{lemma}[theorem]{Lemma}

\newtheorem{proposition}[theorem]{Proposition}
\newtheorem{corollary}[theorem]{Corollary}

\theoremstyle{definition}
\newtheorem{definition}[theorem]{Definition}
\newtheorem{remark}[theorem]{Remark}
\newtheorem{example}[theorem]{Example}
\newtheorem{problem}[theorem]{Problem}

\newcommand{\dda}{\mathord{\mbox{\makebox[0pt][l]{\raisebox{-.4ex}{$\downarrow$}}$\downarrow$}}}

\newcommand{\rom}[1]{\rm{\uppercase\expandafter{\romannumeral #1}}}

\newcommand{\fin}{\mathrm{fin}}

\makeatletter
\def\ps@pprintTitle{%
	\let\@oddhead\@empty
	\let\@evenhead\@empty
	\def\@oddfoot{\reset@font\hfil\thepage\hfil}
	\let\@evenfoot\@oddfoot
}
\makeatother

\begin{document}
	\begin{frontmatter}
		
		\title{Posets Uniquely Determined by its Scott Compact Saturated Subsets \tnoteref{t1}}
		\tnotetext[t1]{This work is supported by the National Natural Science Foundation of China (No.12231007) }
		\author{Huijun Hou}
		\ead{houhuijun2021@163.com}
		\author{Qingguo Li\corref{a1}}
		\address{School of Mathematics, Hunan University, Changsha, Hunan, 410082, China}
		\cortext[a1]{Corresponding author.}
		\ead{liqingguoli@aliyun.com}

	\begin{abstract}
	Inspired by   Zhao and Xu's study on which a dcpo   can be  determined by its Scott closed subsets lattice, we further  investigate whether a poset (or dcpo) $P$ is able to be determined by the family  $\mathcal Q(P)$ of its Scott compact saturated subsets, in the sense  that the isomorphism between $(\mathcal Q(P), \supseteq)$ and  $(\mathcal Q(M), \supseteq)$ implies the isomorphism between  $P$ and $M$ for any poset (or dcpo) $M$, in such case, $P$ is called $\mathcal Q_{\sigma}$-unique.
	Quasicontinuous domains are proved to be $\mathcal Q_{\sigma}$-unique posets and draw support from which, we provide a class of $\mathcal Q_{\sigma}$-unique dcpos. We also define a new kind of posets called $K_D$ and show that every co-sober  $K_D$ poset is $\mathcal Q_{\sigma}$-unique. It  even yields another kind of $\mathcal Q_{\sigma}$-unique dcpos. It is gratifying that  weakly well-filtered co-sober posets are also   $\mathcal Q_{\sigma}$-unique. At last, we distinguish among the conditions which make  a poset (or dcpo) $\mathcal Q_{\sigma}$-unique from each other by some examples; meanwhile, it is confirmed  that none of them except the property of being co-sober are necessary for a poset (or dcpo) to be $\mathcal Q_{\sigma}$-unique.
	
	\end{abstract}
\end{frontmatter}
	
\section{Introduction}
As defined by Drake and Thron in \cite{drake}, the closed subsets lattice $\Gamma(X)$ of  a topological space $X$ represents a lattice $L$ if there is a lattice isomorphism, where the order on $\Gamma(X)$ is the inclusion and  $X$ is called a representation space of  $L$.
Besides, they  proved that  a lattice has  a topological
representation iff it is a $\mathcal C$-lattice.
 We regard the arbitrary two representation spaces of a $\mathcal C$-lattice as the same in the sense of homeomorphisms. All distinct $T_0$ representation spaces of a $\mathcal C$-lattice  compose a family called the representation family. Drake and Thron also reached a crucial result that every representation family of a $\mathcal C$-lattice $L$ have exactly
one element $X$ if and only if every irreducible element of $L$ is strongly irreducible,
that is to say, $X$ has the property that the order isomorphism between $\Gamma (X)$ and  $\Gamma (Y)$ for any $T_0$ space $Y$ implies the  homeomorphism between $X$ and $Y$, which is called $\Gamma$-determined in \cite{zhao2022}, if and only if every irreducible element of $\Gamma(X)$ is strongly irreducible. Further, one deduce from which immediately  that a $T_0$ space $X$ is $\Gamma$-determined if and only if $X$ is sober and $T_D$ (see \cite[Definition 2.1]{thron}), which is first explicitly stated in \cite{td}
(line 11-13, page 504). This indeed identifies a class of $T_0$ spaces which are  determined by their closed subsets lattices.

Inspired by the aforementioned results, Zhao and Xu further focused on the directed complete posets, i.e. dcpos to study whether a dcpo can be uniquely determined by its Scott closed subsets lattice \cite{zhao2018,xls}. They called  $L$ a  $C_{\sigma}$-unique dcpo if for any dcpo $M$, $L$ is isomorphic to $M$ when $\Gamma(L)$ is isomorphic to $\Gamma(M)$. Meanwhile, they provided some sufficient conditions to ensure that a dcpo is $C_{\sigma}$-unique. But unfortunately, none of such conditions are necessary for a dcpo to be $C_{\sigma}$-unique. The full characterization of a $C_{\sigma}$-unique
dcpo is still open.

Hofmann-Mislove's Theorem  is much valuable for demonstrating  the close connections between domain theory and topologies, as worked out in \cite[Section IV-2]{clad}, the open sets lattice $\mathcal O(X)$ and  $\mathcal Q(X)$,  the family of all compact saturated subsets,  are always Lawson duals when $X$ is a locally compact sober space. Additionally, Escard$\mathrm{\acute{o}}$ et al. in \cite{cosober} derived that for each sober and co-sober core compactly generated space $X$, $\mathcal O(X)$ is the Lawson dual of $\mathcal Q(X)$, which provided an alternative condition for the dual Hofmann-Mislove Theorem. It is worth saying that
compact saturated subsets  play a vital role in either way. Meanwhile, considering of  the established facts that each compact saturated subset of a topological space $X$  is closed in the co-compact topology, and that all closed subsets of $X$ could uniquely determine $X$ under some conditions,  it would be meaningful to investigate whether a topological space is able to be uniquely determined by its compact saturated subsets.

In this paper, we focus on the posets (or  dcpos) to explore when the family of their Scott compact saturated subsets are capable of determining them.  We call  $P$ a  $\mathcal Q_{\sigma}$-unique poset (or dcpo) if for any poset (or dcpo) $M$, $P$ is isomorphic to $M$ provided that $\mathcal Q(P)$ is isomorphic to $\mathcal Q(M)$, where $\mathcal Q(P)$ and $\mathcal Q(M)$ are endowed with the reverse inclusion order.  We first show that each quasicontinuous domain is a $\mathcal Q_{\sigma}$-unique poset by drawing the support from the $\mathcal Q$-determinacy introduced by Shen in \cite{sa}. Then based on the quasicontinuity, we give a class of  $\mathcal Q_{\sigma}$-unique dcpos.
We also introduce a  concept  called a $K_D$ poset and   identify a class of $\mathcal Q_{\sigma}$-unique posets which are $K_D$ and co-sober. In addition, the $\mathcal Q$-determinacy and  the co-sobriety  together form a sufficient condition to make a poset $\mathcal Q_{\sigma}$-unique, which immediately suggests that every weakly well-filtered and co-sober poset is $\mathcal Q_{\sigma}$-unique. At last, we utilize some examples to show that the sufficient conditions  making a poset (or dcpo)  $\mathcal Q_{\sigma}$-unique are different from each other and except the  co-sobriety, none of them is necessary for a poset (or dcpo) to be $\mathcal Q_{\sigma}$-unique. The necessity of the co-sobriety is still unclear.

\section{Preliminaries}
We will briefly introduce some basic concepts in domain theory. Readers are  referred to \cite{clad} and \cite{nht2} for more details.

Given a poset $P$. A subset $A$ of $P$ is called \emph{upper set} if $A = {\uparrow}A$, where ${\uparrow}A = \{x\in P: x\geq a \ \mathrm{for\ some\ a\in A}\}$. The \emph{lower set} is defined dually.  In particular, we write ${\uparrow} \{x\}$ as ${\uparrow}x$ and ${\downarrow}\{x\}$ as ${\downarrow}x$, respectively. $D\subseteq P$ is \emph{directed} if every finite subset of $D$ has an upper bound in $D$. We call a subset $U$ of $P$ \emph{Scott open} if  $U = {\uparrow} U$ and  for any directed subset $D$  whose supremum $\sup D$ exists, $\sup D\in U$ implies $D\cap U\neq\emptyset$. All Scott open subsets of $P$ form a topology, which is called \emph{Scott topology} and denoted by $\sigma(P)$.  Let $\Sigma P$ denote the topological space $(P, \sigma(P))$. A mapping $f$ between the posets $P$ and $M$ is called \emph{Scott-continuous}  if $f(\sup D) = \sup f(D)$ holds for each directed subset $D$ of $P$ with $\sup D$  existing, which is indeed a continuous mapping from $\Sigma P$ to $\Sigma{M}$.

Let $L$ be a dcpo and $G, H$ subsets of $L$. Recall that the \emph{Smyth order} on $L$ is defined as $G\sqsubseteq H \Longleftrightarrow H\subseteq {\uparrow}G$.
We say that $G$ is \emph{way below} $H$, it is written as $G \ll H $, if
for every directed subset $D\subseteq L$, $\sup D \in {\uparrow}H $ implies $d\in {\uparrow}G $ for some $d \in D$.
We write $G \ll x$ for $G \ll\{x\}$. Note that $y\ll x $ is
unambiguously defined when $G$ is taken as the singleton $\{y\}$. $L$ is called a \emph{quasicontinuous domain} if for each $x \in L $,
$\fin(x) = \{F: F\ \mathrm{is\ finite}, F \ll x\}$
is a directed family with respect to the Smyth order and ${\uparrow}x = \bigcap\{{\uparrow}F: F\in \fin(x)\}$. $L$ is called a \emph{domain} if the set $\dda x = \{y: y\ll x\}$ is directed and $x = \sup \dda x$.

For a topological space  $X$, the specialization preorder $\leq$ defined by $x \leq y$ iff $x \in cl(\{y\})$
 will be a partial order when $X$ is a $T_0$ space; in such a case, it is called \emph{specialization order}. Clearly, the fact that the Scott topology is $T_0$ guarantees that the  specialization order with respect to the Scott topology is the partial order.
A subset of $X$ is \emph{saturated} if it is  the intersection of its open neighborhoods, or equivalently,  it is an upper
set in the order of specialization.  Let $\mathcal Q(X)$ denote the set of all  compact saturated subsets of $X$.  One can easily see that $\mathcal Q(X)$ is closed under finite unions, thus constitutes a closed basis and the topology generated by which is called  co-compact topology.
 $X$ is called \emph{well-filtered}  if for each filter basis $\mathcal C$ of $\mathcal Q(X)$  and any open subset $U$ of $X$ with $\bigcap \mathcal C\subseteq U$, $K \subseteq U $ for some $K\in \mathcal C$.
We use $\Gamma(X)$ to denote the family of all closed subsets.   $A\in \Gamma(X)$ is called \emph{irreducible} if for any $B, C\in \Gamma(X)$, $A\subseteq B\cup C$ implies that $A\subseteq B$ or $A\subseteq C$. $X$ is called \emph{sober} if every nonempty irreducible closed subset is the closure of a point.	In particular, each sober space is indeed a dcpo under its specialization order.

\begin{definition}\cite{cosober}
\begin{enumerate}
\item A compact saturated subset $K$ of a topological space $X$ is called \emph{$k$-irreducible} if it cannot be written as the union of two proper compact saturated subsets.
    \\
 Obviously, for every $x\in X$, $\uparrow{x}$ is  compact saturated
and $k$-irreducible. We use $kIRR(X)$ to denote the set of all $k$-irreducible subsets of $X$.
\item We call a space $co$-$sober$ if every $k$-irreducible compact saturated subset is  the saturate of a point.
\end{enumerate}
	
\end{definition}

\begin{remark}
Let $X$ be a topological space. 
\begin{enumerate}
\item For each compact saturated subset $K$ of $X$, if $K$ is irreducible closed in the co-compact topology, then it  is $k$-irreducible. 
    \item It is worth noting that there may be subsets of $X$ that are irreducible closed in the co-compact topology but not $k$-irreducible. Thus a co-sober space  may not be a sober space when  it endows with the co-compact topology. We will provide an example below to illustrate this statement.
\end{enumerate}
\end{remark}

\begin{example}
Let $\mathbb N$ denote the set of all natural numbers and $L = \{a, b\}\cup\{-n: n\in \mathbb{N}\}$. For any $x, y\in L$, the order $x\leq y$ is defined as below:
\begin{itemize}
 \item  $x=a, y = -n$, for all $n\in\mathbb N$; or
  \item  $ x=b, y = -n$, for all $n\in\mathbb N$; or
  \item  $x=-n, y=-m$, where $m, n
 \in \mathbb N$ and $m\leq n$.
 \end{itemize}
      \begin{figure}[h]
		\centering
		\begin{tikzpicture}[scale=1]
		\draw (3,1)  circle (3pt) node[left] {-1};
\draw (3,0)  circle (3pt) node[left] {-2};
		\draw (3,-1)  circle (3pt) node[left] {-3};
		\draw (3,-2)  circle (3pt) node[left] {-n};
		
	\draw (3,0.9)--(3,0.1);
	\draw (3, -0.1)--(3,-0.9);
	\draw (3,-1.1)--(3,-1.9);	
\draw[densely dashed] (3,-2.1)--(3,-3);	
	\draw (2,-3.5)  circle (3pt) node[below] {$a$};
\draw (4,-3.5)  circle (3pt) node[below] {$b$};
			
		\end{tikzpicture}
		\scriptsize    \\ \small{Figure 1: The dcpo $ L$ of Example 2.3.}
	\end{figure}
Then $L$ is  a dcpo (see Figure 1). It is easily seen that $\mathcal Q(L) = \{L\}\cup\{\uparrow x: x\in L\}$, which implies that
\begin{center}
$kIRR(L)= \{\uparrow x: x\in L\}$.
\end{center}
Thus, $\Sigma L$ is co-sober.

Now we prove that $L$ is not sober when it is endowed  with the co-compact  topology. Let $\Gamma_{co}(L)$ denote the set of all closed subsets in the co-compact topology of $\Sigma L$. Then we have
\begin{center}
  $\Gamma_{co}(L)=\{L\}\cup\{\uparrow x: x\in L\}\cup\{\uparrow a \;\cap \uparrow b\}$.
\end{center}
 One can easily verify that $\uparrow a \;\cap \uparrow b =\{-n: n\in \mathbb N\}$   is an irreducible closed subset of $L$ in the co-compact topology but not a compact saturated subset of $\Sigma L$, thus naturally, it is not $k$-irreducible. In addition, $\uparrow a \;\cap \uparrow b$ can not be written as the closure of a singleton set, since there is no  greatest element in $\uparrow a \;\cap \uparrow b$ under the order generated by the co-compact topology. Hence, $L$ endowed with  the co-compact topology is not sober.

\end{example}

\begin{definition}\cite{ww}
	A topological space $X$ is called \emph{weakly well-filtered}, whenever a nonempty open set $U$ contains a filtered intersection $\bigcap_{i\in I}Q_{i}$ of compact saturated subsets, then $U$ contains $Q_{i}$ for some $i\in I$.
\end{definition}
	
	\section{Main results}
In this section,   without specifying, the topology on a poset $P$ is always the Scott topology. $\mathcal Q(P)$ denotes the collection composed of all compact saturated subsets.  Conventionally, the order on $\mathcal Q(P)$  is the reverse inclusion.

\begin{definition}
	A poset (or dcpo) $P$ is called \emph{$\mathcal Q_{\sigma}$-unique} if for any poset (or dcpo) $M$, there will be an order isomorphism between $P$ and $M$ provided  that $\mathcal Q(P)$ is  isomorphic to $\mathcal Q(M)$.
\end{definition}

Recall that He and Wang once raised the notions of $T$-lattices, and proved that if  $P$ and $Q$ are $T$-lattices,  $P$ is order isomorphic to $Q$ if and only if  $\mathcal Q(P)$ is order isomorphic to $\mathcal Q(Q)$ \cite{wky}. Now we generalize the notions to  posets.

\begin{definition}
	We shall say $P$ a \emph{$K_D$ poset} if for every $x\in P$, ${\uparrow}x\setminus\{x\}$ is still compact in $\Sigma P$.
\end{definition}

Note that if ${\uparrow}x\setminus\{x\}$ is compact in $\Sigma P$, then it is naturally closed in the co-compact topology  of $\Sigma P$. It indicates that a $K_D$ poset $P$ must be a $T_D$ poset when $P$ is endowed with the co-compact topology.  However, for a poset, even a dcpo,  though it is $T_D$ in the co-compact topology, we can't infer that it is $K_D$. One can refer to the next example:
\begin{example}
Let $L = \{a, b\}\cup\mathbb{N}$ with the order $x\leq y$ defined as:
\begin{itemize}
 \item  $x=a, y = n$, for all $n\in\mathbb N$; or
  \item  $ x=b, y = n$, for all $n\in\mathbb N$.
 \end{itemize}
      \begin{figure}[h]
		\centering
		\begin{tikzpicture}[scale=1]
	\draw (-2,-1)  circle (3pt) node[above] {1};
\draw (0,-1)  circle (3pt) node[above] {2};
\draw (2,-1)  circle (3pt) node[above] {3};
	\draw (4,-1)  circle (3pt) node[above] {n};
	\draw (-0.05,-2.9)--(-2,-1.1);
	\draw (0,-2.9)--(0,-1.1);
\draw (0.05,-2.9)--(1.95,-1.1);
\draw (0.11,-2.95)--(3.92,-1.05);

\draw (1.9,-2.95)--(-1.93,-1.07);
	\draw (1.95,-2.9)--(0.05,-1.1);
\draw (2,-2.9)--(2,-1.1);
\draw (2.06,-2.9)--(4,-1.1);
	\draw[densely dashed] (4.3,-1)--(5.1,-1);	
\draw[densely dashed] (2.3,-1)--(3.65,-1);	
	\draw (0,-3)  circle (3pt) node[below] {$a$};
\draw (2,-3)  circle (3pt) node[below] {$b$};
			
		\end{tikzpicture}
		\scriptsize    \\ \small{Figure 2: The dcpo $ L$ of Example 3.3.}
	\end{figure}

$L$ is clearly a dcpo (see Figure 2). We first verify that $L$ is  $T_D$ when it is endowed  with the co-compact topology of $\Sigma L$. For an arbitrary $n\in \mathbb N$, $\downarrow n\setminus\{n\}$ in the co-compact topology is actually $\uparrow n\setminus\{n\} = \emptyset$, which is  obviously compact in $\Sigma L$. Thus, each $\downarrow n\setminus\{n\}$ is closed in the co-compact topology. Now consider $\downarrow\{a\}\setminus\{a\}$ in the co-compact topology, it actually equals to the set ${\uparrow}a\cap{\uparrow}b$. The closeness of ${\uparrow}a$ and ${\uparrow}b$ indicates the closeness of ${\uparrow}a\cap{\uparrow}b$ in the co-compact topology. Thus, $\downarrow a\setminus\{a\}$ is closed in the co-compact topology. The case of $b$ is the same as that of $a$, so we don't repeat anymore.  Therefore, we  conclude that $L$ is a $T_D$ dcpo under the co-compact topology.

However, $L$  is not a $K_D$ dcpo, since there is $a\in L$, ${\uparrow a}\setminus\{a\}= \mathbb N$ is not compact in $\Sigma L$.

\end{example}

\begin{lemma}{\rm\cite{jia}}\label{jia}
	Let $K$ be a compact subset of a dcpo $P$. Then every element $x\in K$ is above some minimal element of $K$.
\end{lemma}
This lemma immediately reveals  that for any $K\in \mathcal Q(P)$, $K = {\uparrow}\min K$, where $\min K$ is the set of all minimal elements of $K$.

Let $kIRR(P)$ denote the set of all $k$-irreducible subsets of the poset $P$ with respect to the Scott topology. The order on $kIRR(P)$ is  inherited from $\mathcal Q(P)$.  For the arbitrary posets $P$ and $M$,  one can easily verify that $kIRR(P)$ will be  order isomorphic to $kIRR(M)$ when $\mathcal Q(P)$ is order isomorphic to $\mathcal Q(M)$. 
\begin{theorem} \label{kl}
	Let $P$ be a dcpo. If it satisfies the following property:

	 $\mathbf{(KL):}$ For every $K\in kIRR(P)$, there exists a directed family $\{K_i: i\in I\}\subseteq kIRR(P)$, in which each ${\uparrow}_{kIRR(P)}\{K_i\} = \{K\in kIRR(P): K\subseteq K_i\}$ is a chain in $kIRR(P)$, such that $K = \sup_{i\in I}K_i$.
\vspace{1em}
\\
	Then $P$ is a $\mathcal Q_{\sigma}$-unique dcpo.
	\begin{proof}
		Given a dcpo $M$. Assume that $\mathcal Q(P)$ is order isomorphic to $\mathcal Q(M)$.  We first claim that if $C\in kIRR(P)$ satisfying ${\uparrow}_{kIRR(P)}\{C\}$ is a chain in $kIRR(P)$, then $C = {\uparrow}c_0$ for some $c_0\in P$.
		Since ${\uparrow}_{kIRR(P)}\{C\}$ is a chain, $\{{\uparrow}c: c\in C\}\subseteq {\uparrow}_{kIRR(P)}\{C\}$ is also a chain, in other words, $C$ is a chain in $P$. By Lemma \ref{jia},  we have $C = {\uparrow}\min C$,  obviously, the cardinality of $\min C$  is 1. Thus there exists a $c_0\in P$ such that   $C = {\uparrow}c_0$.

		Now consider an arbitrary  $K\in kIRR(P)$, by Property $\mathbf{(KL)}$, it is the supremum of a directed family $\{K_i: i\in I\}\subseteq kIRR(P)$, where each ${\uparrow}_{kIRR(P)}\{K_i\}$ is a chain in $kIRR(P)$. From the above claim, for every $i\in I$, there exists a $k_i\in P$ such that $K_i = {\uparrow}k_i$.
		Then $K = \sup_{i\in I}K_i  = \sup_{i\in I}{\uparrow}k_i = {\uparrow}\sup k_i$, where the last equation holds owing to the directedness of $\{k_i: i\in I\}$. Thus we can conclude that for every $K\in kIRR(P)$, there is a $x_K\in P$ such that $K = {\uparrow}x_K$, that is to say, $kIRR(P) = \{{\uparrow}x: x\in P\}$.
		
		Since $\mathcal Q(P)$ is order isomorphic to $\mathcal Q(M)$,  $kIRR(P)$ is order isomorphic to $kIRR(M)$. Thus $M$ also has Property $\mathbf{(KL)}$ and   $kIRR(M) = \{{\uparrow}m: m\in M\}$.  It  follows that $P\cong kIRR(P)\cong kIRR(M)\cong M$.
	\end{proof}
\end{theorem}

\begin{definition}\cite{sa}
	Let $M$ be a poset. An element $a\in M$ is called \emph{strongly prime} if for any $K\in \mathcal Q(M)$ with $\wedge K$ existing, $a\geq\wedge K$ implies $a\in K$.
	
	Let $SP(M)$ denote the set of all strongly prime elements in $M$.
\end{definition}

\begin{definition}\cite{sa}
	A poset $M$ is called \emph{$\mathcal Q$-determined} if for any $\mathcal{K}\in \mathcal Q(\mathcal Q(M))$ with $\bigwedge \mathcal{K}$ existing, $\bigwedge \mathcal{K} = \bigcup \mathcal{K}$.
\end{definition}

\begin{lemma}\label{sko}{\rm\cite{sa}}
	If $P$ and $M$ are both $\mathcal Q$-determined posets, then $P$ is order isomorphic to $M$ if and only if $\mathcal Q(P)$ is order isomorphic to $\mathcal Q(M)$.
\end{lemma}

\begin{lemma}\label{wd}{\rm\cite{sa}}
	A poset $M$ is $\mathcal Q$-determined if one of the following two alternatives holds:
	\begin{enumerate}
		\item $M$ is weakly well-filtered;
		\item $\mathcal Q(M)$ is a domain.
	\end{enumerate}
\end{lemma}

Based on the aforementioned results, we reach the following conclusion.

\begin{proposition}\label{qua}
	If $L$ is a quasicontinuous domain, then $L$ is a $\mathcal Q_{\sigma}$-unique poset.
	\begin{proof}
		From \cite[Proposition I-1.24.2]{clad}, we have that $\mathcal Q(L)$ is a domain if $L$ is a quasicontinuous domain. Thus for any poset $M$ that satisfies $\mathcal Q(L)\cong \mathcal Q(M)$, $\mathcal Q(M)$ is also a domain. Then by Lemmas \ref{sko}  and \ref{wd}, $L\cong M$.
	\end{proof}
\end{proposition}

\begin{theorem}\label{qdcpo}
	Let $L$ be a dcpo  satisfying the conditions below:
	\begin{enumerate}
		\item $\Sigma L$ is co-sober;
		\item for each $x\in L$, there exists a directed subset $D$, where each element $d\in D$  satisfies that ${\uparrow}d$ is a quasicontinuous dcpo,  such that $x = \sup D$.
	\end{enumerate}
	Then $L$ is a $\mathcal Q_{\sigma}$-unique dcpo.
	\begin{proof}
		Assume $M$ is a poset and there exists an order isomorphism $F: \mathcal Q(L)\rightarrow \mathcal Q(M)$. Let $a$ be an element  with  the property that ${\uparrow}a$ is a quasicontinuous dcpo. Note that ${\uparrow}a$ and $F({\uparrow}a)$ are both posets with respect to the
		restricted  order on $L$ and $M$, respectively. Now  we define $F\mid_{{\uparrow}a}: \mathcal Q({\uparrow}a)\rightarrow \mathcal Q(F({\uparrow}a))$ as:
		\begin{center}
			for any $K\in \mathcal Q({\uparrow}a)$, $F\mid_{{\uparrow}a}(K) = F(K)$.
		\end{center}
		
		Claim 1: $F\mid_{{\uparrow}a}$ is well-defined.
		
		Let $K$ be a compact saturated subset of ${\uparrow}a$. From the fact that the Scott topology  on ${\uparrow}a$ is actually inherited from which of $L$, we have  $K\in\mathcal Q(L)$. It implies $F(K)\in\mathcal Q(M)$. Since $K\subseteq {\uparrow}a$ and $F$ is order preserving, we have $F(K)\subseteq F({\uparrow}a)$. It follows that $F(K)\in\mathcal Q(F({\uparrow}a))$. Thus $F\mid_{{\uparrow}a}$ is well-defined.
		
		Claim 2: $F\mid_{{\uparrow}a}$ is an order embedding.
		
		It is obvious owing to the fact that $F$ is an order embedding.
		
		Claim 3:  $F\mid_{{\uparrow}a}$ is surjective.
		
		Let $K'\in \mathcal Q(F({\uparrow}a))$. Then $K'\in\mathcal Q(M)$. As $F$ is surjective, there exists $K_{0}\in \mathcal Q(L)$ such that $F(K_0) = K'\subseteq F({\uparrow}a)$.   By the property of being an order embedding again, $K_0\subseteq {\uparrow}a$. Thus $K_0\in \mathcal Q({\uparrow}a)$. It entails that $F\mid_{{\uparrow}a}(K_0) = K'$, that is,  $F\mid_{{\uparrow}a}$ is surjective.
		
		So $F\mid_{{\uparrow}a}$ is an order isomorphism. Proposition \ref{qua} tells us such ${\uparrow}a$ is a $\mathcal Q_{\sigma}$-unique poset. Thus $\mathcal Q({\uparrow}a)\cong \mathcal Q(F({\uparrow}a))$ implies ${\uparrow}a\cong F({\uparrow}a)$. Hence, there is a $y_a\in M$ such that $F({\uparrow}a) = {\uparrow}y_a$.
		
		Now let $x$ be an arbitrary element in $L$. By (2),  $x = \sup D$ for some directed subset $D$ of $L$ in which each ${\uparrow}d$ is a quasicontinuous  dcpo. Based on the above, there exists a $y_d\in M$ such that $F({\uparrow}d) = {\uparrow}y_d$ for each $d\in D$. Then the following is derived:
		\begin{align*}
			F({\uparrow}x) &= F({\uparrow}\sup D) \\ &= F(\sup\! _{Q(L)} \{{\uparrow}d: d\in D\} ) \\ & =  \sup\!_{Q(M)} \{F({\uparrow}d): d\in D\} \\ &= \sup\!_{Q(M)}\{{\uparrow}y_{d}: d\in D \}  \\ & =  {\uparrow}\sup\!_{d\in D} y_{d}.
		\end{align*}
		where the last equation holds due to the facts that $\{y_d: d\in D\}$ is  directed  and  $M$ is a dcpo. Besides, we write $\sup\!_{d\in D} y_{d}$ as $y_x$. Then we can define
		a map $f: L\rightarrow M$ as
		\begin{center}
			for any $x\in L$, $f(x) = y_x$,
		\end{center}
		 which clearly satisfies   ${\uparrow}f(x) = F({\uparrow}x)$. Besides, it is easy to verify that $f$ is well-defined and  an order embedding. The remainder that needs to be proved is the surjectivity of $f$.
		
		For each $m\in M$, clearly, ${\uparrow}m$ is a $k$-irreducible subset of $M$. By the fact that $F$ is an order isomorphism, we know $F\mid_{kIRR(L)}: kIRR(L)\rightarrow kIRR(M)$ is also an isomorphism. Thus there must exist a $C\in kIRR(L)$ such that $F(C) = {\uparrow}m$. Because $\Sigma L$ is co-sober, $C = {\uparrow}c$ for some $c\in L$. It means $F({\uparrow}c) = {\uparrow}m$. So ${\uparrow}f(c) = {\uparrow}m$; hence, $f(c) = m$. Thus $f$ is surjective.
		
	\end{proof}
\end{theorem}

Note that for an element $x$ in a dcpo $L$, if ${\uparrow}x$ is a chain, then it must be a continuous lattice and thus a quasicontinuous dcpo.  So naturally, the following conclusion emerges:
\begin{corollary}\label{chain}
	Let $L$ be a dcpo that satisfies the conditions as follows:
	\begin{enumerate}
		\item $\Sigma L$ is co-sober;
		\item for each $x\in L$, there exists a directed subset $D$, where  ${\uparrow}d$ is a chain in $P$ for each  $d\in D$,  such that $x = \sup D$.
	\end{enumerate}
	Then $L$ is a $\mathcal Q_{\sigma}$-unique dcpo.
\end{corollary}
It is worth noting that Condition (2) of Corollary \ref{chain} and Property $\mathbf{(KL)}$ in Theorem \ref{kl} are inextricably related to the chain, which arouses our interest to study whether there is a possibility that they are the same.

\begin{lemma}\label{kl2}
	If a dcpo $L$ possesses Property $\mathbf{(KL)}$, then it  satisfies Condition (2) of Corollary \ref{chain}.
	\begin{proof}
		Let $x$ be an arbitrary element in $L$. Then ${\uparrow}x\in kIRR(L)$. By Property $\mathbf{(KL)}$, there exists a directed family $\mathcal K$ contained in $kIRR(L)$, in which  each ${\uparrow}_{kIRR(L)}\{K\}$ is a chain in $kIRR(L)$, such that ${\uparrow}x = \sup\mathcal K$. As the claim of Theorem \ref{kl} shows,  each $K\in \mathcal K$ is actually a chain and there exists a $x_K\in L$ such that $K = {\uparrow}x_K$.
		The directedness of $\mathcal K$ implies that the set $\{x_K: K\in \mathcal K\}$ is  directed.  Thus $\sup_{K\in \mathcal K}x_K$ exists and ${\uparrow}x = \sup\{{\uparrow}x_K: K\in \mathcal K\} = {\uparrow}\sup_{K\in \mathcal K}x_K$, it is equivalent to saying that $x = \sup_{K\in \mathcal K}x_K$. So Condition (2) of Corollary \ref{chain} is satisfied.
	\end{proof}
\end{lemma}

\begin{proposition}
	Let $L$ be a dcpo. Then $L$ has Property $\mathbf{(KL)}$ if and only if it satisfies the conditions of Corollary \ref{chain}.
	\begin{proof}
		$(\Rightarrow)$: Lemma \ref{kl2} has shown that $L$ satisfies Condition (2) of  Corollary \ref{chain}. See the proof of Theorem \ref{kl}, we find that $kIRR(L) = \{{\uparrow}x: x\in L\}$ whenever $L$ has Property $\mathbf{(KL)}$. It means that $\Sigma L$ is co-sober.
		
		$(\Leftarrow)$: Since $\Sigma L$ is co-sober, $kIRR(L) = \{{\uparrow}x: x\in L\}$.  For each $x\in L$, by Condition (2) in Corollary \ref{chain}, there exists a directed subset $D\subseteq L$ in which each ${\uparrow}d$ is a chain in $L$ such that $x = \sup D$. It implies that ${\uparrow}x = {\uparrow}\sup D = \sup_{kIRR(L)}\{{\uparrow}d: d\in D\}$. By the fact  that each ${\uparrow}d$ is a chain, we have that ${\uparrow}_{kIRR(L)}\{{\uparrow}d\} = \{{\uparrow}c: c\in {\uparrow}d\}$ is a chain in $kIRR(L)$. As a result, $L$ has Property $\mathbf{(KL)}$.
	\end{proof}
\end{proposition}

\begin{theorem}\label{qdcpo2}
	Let $L$ be a dcpo satisfying the following two conditions:
	\begin{enumerate}
		\item $\Sigma L$ is co-sober;
		\item for each $x\in L$, there exists a directed subset $D$ of $L$ such that $x = \sup D$, where each $d\in D$ has the property that ${\uparrow}d\setminus \{d\}$ is compact.
	\end{enumerate}
	Then $L$ is a $\mathcal Q_{\sigma}$-unique dcpo.
	\begin{proof}
		Let $M$ be a dcpo. Assume that there exists an order isomorphism $F$ from $\mathcal Q(L)$ to $\mathcal Q(M)$. Consider $a\in L$ that ${\uparrow}a\setminus \{a\}$ is still compact.  Clearly, ${\uparrow}a\setminus \{a\}\in\mathcal Q(L)$ and  is properly contained in ${\uparrow}a$. Then by the fact that $F$ is an order embedding,  we have  that $F({\uparrow}a\setminus \{a\})$ is properly contained in $F({\uparrow}a)$. Thus there must exist a $y_a$ which belongs to $F({\uparrow}a)$ but not to $F({\uparrow}a\setminus \{a\})$. Then we claim that $F({\uparrow}a) = {\uparrow}y_a$. Obviously, $y_a\in F({\uparrow}a)$ hints that ${\uparrow}y_a\subseteq F({\uparrow}a)$, which is equivalent to saying that $F^{-1}({\uparrow}y_a)\subseteq {\uparrow}a$. Conversely, we just need to show ${\uparrow}a \subseteq F^{-1}({\uparrow}y_a)$. Suppose that $a\notin F^{-1}({\uparrow}y_a)$. Then $F^{-1}({\uparrow}y_a)\subseteq {\uparrow}a\setminus \{a\}$, which implies $y_a\in F({\uparrow}a\setminus \{a\})$, a contradiction. So  ${\uparrow}a \subseteq F^{-1}({\uparrow}y_a)$ and then $F({\uparrow}a) = {\uparrow}y_a$ holds.
		
		For each $x\in L$, by Condition (2), there exists a subset $D\subseteq L$ directed in which each ${\uparrow}d\setminus \{d\}$ is compact such that $x = \sup D$. Assume that  $F({\uparrow}d) = {\uparrow}y_d$ for each $d\in D$. Then we have
		\begin{align*}
			F({\uparrow}x) &= F({\uparrow}\sup D) \\ &= F(\sup\! _{Q(L)} \{{\uparrow}d: d\in D\} ) \\ & =  \sup\!_{Q(M)} \{F({\uparrow}d): d\in D\} \\ &= \sup\!_{Q(M)}\{{\uparrow}y_{d}: d\in D \}  \\ & =  {\uparrow}\sup\!_{d\in D} y_{d},
		\end{align*}
		where the directedness of $\{y_d: d\in D\}$ and the fact that $M$ is a dcpo together guarantee that the last equation holds. Now we  define $f: L\rightarrow M$ such that $f(x) = y_x$ iff ${\uparrow}y_x = F({\uparrow}x)$. One can easily verify that $f$ is an order embedding. The surjection of $f$ can be proved similar to that in Theorem \ref{qdcpo}.
	\end{proof}
\end{theorem}

Given that we have provided  three sufficient conditions  to make a dcpo $\mathcal Q_{\sigma}$-unique (Theorem \ref{qdcpo}, \ref{qdcpo2} and Corollary \ref{chain}). It is therefore necessary to illustrate that the three conditions are different from each other and consider whether there is any property among them being necessary for a dcpo to be  $\mathcal Q_{\sigma}$-unique. Concretely, see examples as follows:
\begin{example}\label{ce}
	(1): Let $L = \mathbb N\cup \{\perp\}$ with the order $x\leq y$ defined as $x = \perp, y\in \mathbb N$ or $x= y \in \mathbb N$ (see Figure 3).
	\begin{figure}[h]
		\centering
		\begin{tikzpicture}[scale=1]
			\draw (0,0)  circle (3pt) node[above] {1};
			\draw (1.5,0)  circle (3pt) node[above] {2};
			\draw (3,0)  circle (3pt) node[above] {3};
			\draw (4.6,0)  circle (3pt) node[above] {n};
			\draw[densely dashed] (3.4,0)--(4.2,0);
			\draw[densely dashed] (5, 0)--(6,0);
			
			\draw (2.28,-1.5)  circle (3pt) node[right] {$\bot$};
			\draw (0,-0.1)--(2.21,-1.42) (1.5,-0.1)-- (2.25,-1.4) (3,-0.1)--(2.3,-1.4) (4.6,-0.1)--(2.35,-1.42);

		\end{tikzpicture}
		\scriptsize    \\ \small{Figure 3: The dcpo $ L$ of (1) in Example 3.16.}
	\end{figure}
	One can easily see that $L$ is a dcpo and $\mathcal Q(L) = \{L\}\cup\{F: F\subseteq \mathbb N\ \textrm{is finite}\}$. Further, $kIRR(L)=\{L\}\cup\{\{n\}: n\in \mathbb N\}$,  which immediately entails that $\Sigma L$ is co-sober and for each $x\in L$,  ${\uparrow}x$   is a quasicontinuous dcpo. Thus it satisfies the conditions of Theorem \ref{qdcpo}.  Note that there is only a directed subset $\{\bot\}$ whose supremum is $\bot$. However, ${\uparrow}\bot = L$ is not a chain and ${\uparrow}\perp\setminus\{\perp\} = \mathbb N$ is not compact in $\Sigma L$. So $L$ fails to meet  Condition (2) of Corollary \ref{chain} and that of Theorem \ref{qdcpo2}, which hints that neither Condition (2) of Corollary \ref{chain}  nor that of  Theorem \ref{qdcpo2} is necessary for a dcpo to be $\mathcal Q_\sigma$-unique.
	
	(2): Let $M = (\mathbb N\cup \{\infty\}) \cup(\mathbb N\times \mathbb N\cup \{\top\})$. For all $x, y\in M$, we define $x\leq y$ if and only if  one of the following
	alternatives holds:
	\begin{itemize}
		\item $ x\leq y$ in $\mathbb N$;
		\item $x\in \mathbb N$, $y = \infty$;
		\item $x= \infty$, $y\in \mathbb N\times \mathbb N\cup \{\top\}$;
		\item $x= (m_1, n_1)$, $y = (m_2, n_2)\in \mathbb N\times \mathbb N$, $m_1 = m_2$, $n_1\leq n_2$;
		\item $x\in \mathbb N\times \mathbb N$,  $y = \top$.
	\end{itemize}
	
	\begin{figure}[h]
		\centering
		\begin{tikzpicture}[scale=1]
			\draw (0,6.3)  circle (3pt) node[right] {\; \Large{$\top$}};
			\draw (-2.4,4)  circle (3pt);
			\fill [black] (-2.4,2.8)  circle (3pt);
			\draw (-2.4,1.6)  circle (3pt) node[left]{(1,1)};
			\draw (-2.4,1.7)--(-2.4,2.7) (-2.4,2.9)--(-2.4,3.9);
			\draw (-0.8,4)  circle (3pt) node[left]{(2,3)};
			\draw (-0.8,2.8)  circle (3pt);
			\fill[black] (-0.8,1.6)  circle (3pt);
			\draw (-0.8,1.7)--(-0.8,2.7) (-0.8,2.9)--(-0.8,3.9);
			\draw (0.8,4)  circle (3pt);
			\draw (0.8,2.8)  circle (3pt);
			\draw (0.8,1.6)  circle (3pt);
			\draw (0.8,1.7)--(0.8,2.7) (0.8,2.9)--(0.8,3.9);
			\fill [black] (2.4,4)  circle (3pt) node[right] at (2.4,3.8){(n,3)};
			\draw (2.4,2.8)  circle (3pt);
			\draw (2.4,1.6)  circle (3pt) ;
			\draw (2.4,1.7)--(2.4,2.7)  (2.4,2.9)-- (2.4,3.9);
			\draw [densely dashed] (1.1,4)--(2.1,4);
			\draw [densely dashed] (1.1,2.8)--(2.1,2.8);
			\draw [densely dashed] (1.1,1.6)--(2.1,1.6);
			\draw [densely dashed] (2.6,4)--(3.75,4);
			\draw [densely dashed] (2.6,2.8)--(3.75,2.8);
			\draw [densely dashed] (2.6,1.6)--(3.75,1.6);
\draw [densely dashed] (-2.4,4.1)--(-2.4,5.1) (-0.8,4.1)--(-0.8,5.1) (0.8,4.1)--(0.8,5.1) (2.4,4.1)--(2.4,5.1);
			\draw [densely dashed] (-2.4,5.1)--(-0.1,6.3) (-0.8,5.1)--(-0.05,6.25) (0.8,5.1)--(0.05,6.22) (2.4,5.1)--(0.1,6.3);
			
			\draw (0,0)  circle (3pt) node[right] {\ \Large{$\infty$}};
			\draw (0,-1.2) circle (3pt)node[right]{\ \large 3};
			\draw (0,-2.2)  circle (3pt)node[right]{\ \large 2};
			\draw(0,-3.2)  circle (3pt) node[right] {\ \large 1};
			\draw (0,-3.1)--(0,-2.3) (0,-2.1)--(0,-01.3);
			\draw [densely dashed] (0,-1)--(0,-0.2);
			\draw  (-0.05,0.1)--(-2.4,1.5) (-0.02,0.1)--(-0.8, 1.5) (0.02,0.1)--(0.8,1.5) (0.05,0.1)--(2.4,1.5);
			
			\draw  [color=red, rounded corners =10pt]
			(-3.4,5.2)--(-3.4,1.25)--(-1.3,1.25)--(-1.35,1.9)--(-0.1,1.9)--(-0.1,1.2)--(1.4,1.2)--(1.4,4.6)--(3.7,4.6)--(3.7,1);	
			\draw  [color=blue, rounded corners =10pt]%
			(-3.1,5.2)--(-3.1,3.1)--(-1.7,3.1)--(-1.7,1)--(1.65,1)--(1.65,4.3)--(3.4,4.3)--(3.4,1);
			\draw  [color=darkgray, rounded corners =10pt]%
			(-2.8,5.2)--(-2.8,3.3)--(-1.5,3.3)--(-1.5,2.2)--(0.2,2.2)--(0.2,1.4)--(4,1.4);
		\end{tikzpicture}
		\scriptsize   \vspace{0.2cm} \\ {The above region of red line is $U_1$, the above region of blue line is $U_2$, the above region of gray line is $U_n$.}
		\\  { Figure 4: The dcpo $M$ of (2) in Example 3.16}.
	\end{figure}
It is easily  seen that $M$  with the above order is a dcpo. For any $x\in \mathbb N\cup (\mathbb N\times \mathbb N\cup \{\top\})$, ${\uparrow}x\setminus\{x\}$ is clearly compact. Focus on the point $\infty$, $\infty = \sup\mathbb N$, where $\mathbb N$ is directed and each $n\in \mathbb N$ satisfies that ${\uparrow}n\setminus\{n\} = {\uparrow}\{n + 1\}$ is compact. Thus the conditions of Theorem \ref{qdcpo2} will be met so long as we show that $\Sigma M$ is co-sober. To do this, we declare that each  compact saturated subset  of $M$ has finite minimal elements. We just need to consider $K\in \mathcal Q(M)$  that $\min K\subseteq \mathbb N\times\mathbb N$. Suppose that there exists $\min K\subseteq \mathbb N\times\mathbb N$  infinite. Let $\mathbb N_K = \{n\in \mathbb N: \exists\ m\in \mathbb N\ \textrm{s.t.}\ (n, m)\in K\}$ and $m_n = \min\{m\in \mathbb N: (n, m)\in K\}$ for each $n\in \mathbb N_K$.
Now  we set $U_n = M\setminus\bigcup_{k\in \mathbb N_K\setminus\{n\}}{\downarrow}(k, m_k)$ for every $n\in \mathbb N_K$ (we mark three in Figure 4 for a better understanding). It is obvious that each $\bigcup_{k\in \mathbb N_K\setminus\{n\}}{\downarrow}(k, m_k)$ is Scott closed, which implies the openness of each $U_n$. Moreover, $K\subseteq \bigcup_{n\in \mathbb N_K} U_n$. By the compactness, there must exist a finite subset $F\subseteq \mathbb N_K$ such that $K\subseteq \bigcup_{n\in F} U_n$. Since $\min K$ is infinite, there are a $n_0$ being strictly larger than every $n$ in $F$ and a $m_{n_0}\in \mathbb N$ that satisfy $(n_0, m_{n_0})\in U_{n_1}$ for some $n_1\in F$. It follows that $(n_0, m_{n_0})$ does not belong to ${\downarrow}(k, m_k)$ for each $k\in \mathbb N_K\setminus \{n_1\}$,  not to ${\downarrow}(n_0, m_{n_0})$, either, a contradiction.
So $\min K$ is finite for all $K\in\mathcal Q(M)$, which immediately results in the co-sobriety of $M$.

Now focus on  the point $\infty$ again. Note that the directed subsets of $M$ which make $\infty$ supremum are only $\{\infty\}$ and $\mathbb N$. However, for every $x\in \mathbb N\cup \{\infty\}$,  ${\uparrow}x\setminus\{x\}$ is not quasicontinuous since there is no finite subset of $M$ except the empty set  way below $\{\top\}$.  Hence, $M$ does not possess the property (2) of Theorem \ref{qdcpo}. So Condition (2) of Theorem \ref{qdcpo} is not necessary to make a dcpo $\mathcal Q_\sigma$-unique.

In fact, the above two examples are enough to reveal that the conditions of Theorem \ref{qdcpo} and  Theorem \ref{qdcpo2} are different from each other and the conditions of Corollary \ref{chain} cannot be derived by that of Theorem \ref{qdcpo} or  Theorem \ref{qdcpo2}.  To illustrate it as clearly as possible that the conditions of   Theorem \ref{qdcpo}, Corollary \ref{chain} and Theorem \ref{qdcpo2} are different from each other, we  need  another example.

(3) Let $N = (0, 1]$ endow with the Scott topology. Clearly, it is a co-sober dcpo and for each $x\in N$, ${\uparrow}x$ is a chain. Thus $N$ satisfies the conditions of Corollary \ref{chain}. However, for each $y\in (0, 1)$,  ${\uparrow}y\setminus\{y\} = (y, 1]$ is not compact. It suggests that there is no element  except 1 in $N$ can be written as a supremum of some directed subset $D$ in which each ${\uparrow}d\setminus\{d\}$ is compact. So $N$ has nothing to do with Condition (2) of Theorem \ref{qdcpo2}, which means that  the conditions of  Theorem \ref{qdcpo2}     cannot be derived by that of Corollary \ref{chain}.

\end{example}

So far, we can say that the conditions of Theorem \ref{qdcpo}, Corollary \ref{chain} and Theorem \ref{qdcpo2} are quitely different from each other.

Next we will be devoted  to investigating when a poset is $\mathcal Q_{\sigma}$-unique.
\begin{theorem}\label{kdp}
Let $P$ be a $K_D$  and  co-sober poset. Then $P$ is a $\mathcal Q_{\sigma}$-unique poset.
\begin{proof}
	Let $M$ be a poset and assume that there is an order isomorphism $F: \mathcal Q(P)\rightarrow \mathcal Q(M)$. Then we can adopt the same method as that in the proof of Theorem \ref{qdcpo2} to prove that for every $x\in P$, there is a $y_x\in M$ such that $F({\uparrow}x) = {\uparrow}y_x$. Define $f: P\rightarrow M$ as $f(x) = y_x$. The rest is  similar to that in the proof of Theorem \ref{qdcpo}.
\end{proof}
\end{theorem}

\begin{lemma}\label{spri}
Let $P$ be a poset. If $K$ is a strongly prime compact saturated subset of $P$, then $K$ must be a saturation of a point.
\begin{proof}
	We first show that $\bigcup \mathcal U\in \sigma(P)$ whenever $\mathcal U\in \sigma(\mathcal Q(P))$. Obviously, $\bigcup \mathcal U$ is an upper set. Let $D$ be a directed subset for which $\sup D$ exists and $\sup D\in \bigcup \mathcal U$. Then there exists $K_0\in \mathcal U$ such that $\sup D\in K_0$, which implies that ${\uparrow}\sup D\subseteq K_0$. Since $\mathcal U$ is an upper set, ${\uparrow}\sup D\in \mathcal U$. It amounts to saying that  $\sup_{d\in D}{\uparrow}d = {\uparrow}\sup D\in \mathcal U$. Note that $\{{\uparrow}d: d\in D\}$ is actually a directed subset of $\mathcal Q(P)$; thus there must exist a $d_0\in D$ such that ${\uparrow}d_0\in \mathcal U$ by the Scott openness of $\mathcal U$. So $d_0\in \bigcup \mathcal U$.
	
	Assume that $K$ is a strongly prime compact saturated subset of $P$. Let $\mathcal K = {\uparrow}\{{\uparrow}k: k\in K\}$. We prove that $\mathcal K$ is a compact saturated subset of $\mathcal Q(P)$ with respect to the Scott topology. The saturation is clear.  Suppose $\mathcal K\subseteq \bigcup_{i\in I}\mathcal U_i$, where $\{\mathcal U_i: i\in I\}$ is a directed family of $\sigma(\mathcal Q(P))$. Then for each $k\in K$, ${\uparrow}k\in \bigcup_{i\in I}\mathcal U_i$, which implies $k\in \bigcup\mathcal U_i$ for some $i\in I$. Thus $K\subseteq\bigcup_{i\in I}\bigcup\mathcal U_i$. As argued above, every $\bigcup\mathcal U_i$ is Scott open, and besides, $\{\bigcup \mathcal U_i: i\in I\}$ is also directed, so    by the compactness of $K$, we have $K\subseteq \bigcup \mathcal U_{i_0}$ for some $i_0\in I$. It follows that $\mathcal K\subseteq \mathcal  U_{i_0}$, that is to say, $\mathcal K$ is compact. It is easily seen that $K = \bigcup \mathcal K = \bigwedge\mathcal K$. Then the property of being strongly prime implies that $K\in \mathcal K$, that is, $K\subseteq {\uparrow}k_0$ for some $k_0\in K$. As a result, $K = {\uparrow}k_0$.

\end{proof}
\end{lemma}

\begin{theorem}
Let $P$ be a $\mathcal Q$-determined and co-sober poset. Then $P$ is a $\mathcal Q_{\sigma}$-unique poset.
\begin{proof}
	Assume that $M$ is a poset and  $\mathcal Q(P)$ is  isomorphic to $\mathcal Q(M)$, i.e., there is an order isomorphism $F: \mathcal Q(P)\rightarrow \mathcal Q(M)$. We first show that $SP(\mathcal Q(P)) = \{{\uparrow}x: x\in P\}$. The inclusion is trivial by Lemma \ref{spri}. For any $x\in P$, let $\mathcal K$ be a compact saturated subset of $\Sigma (\mathcal Q(P))$ for which  $\bigwedge \mathcal K$ exists and $\bigwedge \mathcal K\leq{\uparrow}x$. Since $P$ is $\mathcal Q$-determined, $\bigwedge \mathcal K = \bigcup\mathcal K$. Then ${\uparrow}x\subseteq \bigcup\mathcal K$; hence, ${\uparrow}x\in \mathcal K$. So $SP(\mathcal Q(P)) = \{{\uparrow}x: x\in P\}$.
	
	One can easily verify that  $F$  restricts to an order isomorphism $F\mid_{SP(\mathcal Q(P))}: SP(\mathcal Q(P))\rightarrow SP(\mathcal Q(M))$. For any  $x\in P$, as argued above, ${\uparrow}x\in SP(\mathcal Q(P))$, which implies that $F({\uparrow}x)\in SP(\mathcal Q(M))$. Since Lemma \ref{spri} reveals that  $SP(\mathcal Q(M))\subseteq \{{\uparrow}m: m\in M\}$, there exists a $m_x\in M$ such that $F({\uparrow}x) = {\uparrow}m_x$ for each $x\in P$. Thus we can define $f: P\rightarrow M$ by $f(x) = m_x$, it is equivalent to saying that ${\uparrow}f(x) = F({\uparrow}x)$. Then  $f$ is an order embedding since $F$ is. Similar to that in Theorem \ref{qdcpo}, we have that $f$ is surjective by the co-sobriety. So $P$ is order isomorphic to $M$, that is, $P$ is a $\mathcal Q_{\sigma}$-unique poset.
	
\end{proof}
\end{theorem}

Based on the above result and with the aid of Lemma \ref{wd}, we derive the following result.
\begin{corollary}\label{wc}
Let $P$ be a weakly well-filtered and co-sober poset. Then $P$ is a $\mathcal Q_{\sigma}$-unique poset.
\end{corollary}

To sum up, the conditions of Lemma \ref{qua}, Theorem \ref{kdp} and Corollary \ref{wc} are all proved to be sufficient  to make a poset   $\mathcal Q_{\sigma}$-unique. What we want to do is to find that the  conditions are necessary for a poset to be  $\mathcal Q_{\sigma}$-unique. Unfortunately, the properties of being  quasicontinuous, $K_D$ and weakly well-filtered fail to be so. Let's make explanations with the following three examples,  respectively.
\begin{example}
(1): See the dcpo $L = \mathbb N\cup \{\bot\}$ given in (1) of Example \ref{ce}. It is obviously a quasicontinuous domain, but  not a  $K_D$ poset since ${\uparrow}\bot\setminus\{\bot\} = \mathbb N$ is not compact in $\Sigma L$. Thus the property of being $K_D$ is not necessary for a poset  to be $\mathcal Q_{\sigma}$-unique.

(2): Let $\mathbb J = \mathbb N\times (\mathbb N \cup \{\omega\})$ be the non-sober dcpo  constructed by Johnstone\cite{john}. 	For any $(j, k), (m, n)\in \mathbb J$, $(j, k)\leq (m, n)$ if and only if one of  the following two alternatives holds (one can see Figure 5 for a better understanding):
\begin{itemize}
	\item $j = m$, $k\leq n$;
	\item $n = \omega$, $k\leq m$.	
\end{itemize}

\begin{figure}
	\centering
	\begin{tikzpicture}[scale=0.8]
		
		\draw (0,0)  circle (3pt)   node[left] {(1,1)} coordinate (a);
		\draw (0,1.9)  circle (3pt)   node[left] {(1,2)} coordinate (b);
		\draw (0,3.8)  circle (3pt)   node[left] {(1,3)} coordinate (c);
		\draw (0,6.5)  circle (3pt)   node[left] {(1,$\omega$)};
		\draw   (0,0.1) -- (0,1.8);
		\draw   (0,2) -- (0,3.7);
		\draw[densely dashed]   (0,3.9) -- (0,6.4);

		\draw (2.5,0)  circle (3pt)   node[left] {(2,1)} coordinate (a);
		\draw (2.5,1.9)  circle (3pt)   node[left] {(2,2)} coordinate (b);
		\draw (2.5,3.8)  circle (3pt)   node[left] {(2,3)} coordinate (c);
		\draw (2.5,6.5)  circle (3pt)   node[left] {(2,$\omega$)};
		
		\draw   (2.5,0.1) -- (2.5,1.8);
		\draw   (2.5,2) -- (2.5,3.7);
		\draw   [densely dashed]   (2.5,3.9) -- (2.5,6.4);

		\draw (5,0)  circle (3pt)   node[left] {\small(3,1)} coordinate (a);
		\draw (5,1.9)  circle (3pt)   node[left] {(3,2)} coordinate (b);
		\draw (5,3.8)  circle (3pt)   node[left] {(3,3)} coordinate (c);
		\draw (5,6.5)  circle (3pt)   node[left] {(3,$\omega$)};
		
		\draw   (5,0.1) -- (5,1.8);
		\draw   (5,2) -- (5,3.7);
		\draw[densely dashed]   (5,3.9) -- (5,6.4);

		\draw (8,0)  circle (3pt)   node[left] {(m,1)} coordinate (a);
		\draw (8,1.9)  circle (3pt)   node[left] {(m,2)} coordinate (b);
		\draw (8,3.8)  circle (3pt)   node[left] {(m,3)} coordinate (c);
		\draw (8,6.5)  circle (3pt)   node[left] {(m,$\omega$)};

		\draw   (8,0.1) -- (8,1.8);
		\draw   (8,2) -- (8,3.7);
		\draw[densely dashed]   (8,3.9) -- (8,6.4);
		\draw[densely dashed] (5.4,0)--(6.6,0) (5.4,1.9)--(6.6,1.9)  (5.4,3.8)--(6.6,3.8)  (5.4,6.5)--(6.6,6.5);
		\draw[densely dashed] (8.4,0)--(10,0) (8.4,1.9)--(10,1.9) (8.4,3.8)--(10,3.8)  (8.4,6.5)--(10,6.5) ;
		
		\draw (-0.8,0.4)--(10,0.4)   (-0.8,-0.4)--(10,-0.4);
		\draw (-0.8,2.3)--(10,2.3)    (-0.8,1.5)--(10,1.5);
		\draw (-0.8,4.2)--(10,4.2)   (-0.8,3.4)--(10,3.4);
		
		\draw (-0.8,0.39) arc (90:270:0.4);
		\draw (-0.8,2.3) arc (90:270:0.4);
		\draw (-0.8,4.2) arc (90:270:0.4);

		\draw  (0.03,6.42)--(0.8,0.42)  ;
		\draw (2.53,6.42)--(3.3,2.26)  ;
		\draw (5.03,6.42)--(5.8,4.2)  ;
		
	\end{tikzpicture}
	\scriptsize   \\
	\vspace{0.2cm} { Figure 5: Johnstone's non-sober dcpo $\mathbb J$, which is in Example 3.21}.
\end{figure}
It has been shown that $\mathcal Q(\mathbb J) = \mathbf 2^{\max \mathbb J}\cup \{{\uparrow}F: F\subseteq \mathbb J{\setminus} \max \mathbb J\ \mathrm{is\ finite}\}$ (see \cite[Example 3.1]{ww}), where $\max \mathbb J = \{(n, \omega): n\in \mathbb N\}$ and $\mathbf 2^{\max \mathbb J}$ is the collection of  all subsets of $\max \mathbb J$. From which we can derive that $kIRR(\mathbb J) = \{{\uparrow}(m, n): (m, n)\in\mathbb J\}$, it indicates the co-sobriety of $\Sigma \mathbb J$. Moreover,
\cite[Example 3.1]{ww} also showed that $\mathbb J$ is weakly well-filtered but not well-filtered. So  $\mathbb J$ is a $\mathcal Q_{\sigma}$-unique poset by Theorem \ref{wc} but not quasicontinuous. As a result, quasicontinuity is not the necessary condition for a poset to be $\mathcal Q_{\sigma}$-unique.

(3): Let $\mathbb S = \mathbb J\cup \{\top\}$. For any $x, y\in \mathbb S$, define $x\leq y$ if and only if $x, y\in \mathbb J, x\leq y$ or $x\in \mathbb J, y = \top$.
In order to show that $\mathbb S$ is $\mathcal Q_{\sigma}$-unique,
we need to distinguish the following three cases:

Case 1: $(m, n)\in \mathbb J\setminus \max \mathbb J$: Then ${\uparrow}(m, n)\setminus\{(m, n)\} = {\uparrow}(m, n+1)$, which is compact in $\Sigma \mathbb S$.

Case 2: $(m, n)\in \max \mathbb J$: Then ${\uparrow}(m, n)\setminus\{(m, n)\} = \{\top\}$, which is compact in $\Sigma \mathbb S$.

Case 3: Consider the element $\top$, ${\uparrow}\top\setminus\{\top\} = \emptyset$ still belongs to $\mathcal Q(\mathbb S)$.

So $\mathbb S$ is a $K_D$ poset. Besides,  similar to that in $\mathbb J$, it is not hard to see that $\mathcal Q(\mathbb S) = \{A\cup \{\top\}: A\in \mathbf 2^{\max \mathbb J}\}  \cup \{{\uparrow}F: F\subseteq \mathbb J\setminus \max \mathbb J\ \mathrm{is\ finite}\}$.
It still means that $\Sigma \mathbb S$ is co-sober. Thus $\mathbb S$ is $\mathcal Q_{\sigma}$-unique by Theorem \ref{kdp}. Now let $K_F = {\uparrow}(\max\mathbb J\setminus F)$ for each finite subset $F$ of $\max\mathbb J$. Then $\{K_F: F\subseteq \max\mathbb J\ \textrm{is finite}\}$ is filtered and $\bigcap_F K_F = \{\top\}\in \sigma(\mathbb S)$. But  no $K_F$  is contained in $\{\top\}$. Thus $\Sigma\mathbb S$ is not weakly well-filtered, which reveals that the weak well-filteredness is not the necessary property to make a poset $\mathcal Q_{\sigma}$-unique.
\end{example}
In particular, $L= \mathbb N\cup \{\bot\}$ given in (1) is also a weakly well-filtered and co-sober poset. So it illustrates that the property  of being $K_D$  cannot be derived by    the weak well-filteredness and co-sobriety. Thus the conditions of Theorem \ref{kdp} and Corollary \ref{wc} are different from each other. In addition, one can prove that $\mathbb J$ is $K_D$ in a similar way to that in (3). So combining with (1), we have that the conditions of Proposition \ref{qua} and those of Theorem \ref{kdp} are different from each other. The aforementioned Johnstone's example $\mathbb J$ in (2)  has shown that the condition of Proposition \ref{qua} cannot be deduced from that of Corollary \ref{wc}. But it remains unclear as to whether the converse deduction is valid, which  leads to the following question for a deeper exploration:

\begin{problem}
Is the  co-sobriety necessary for a poset to be $\mathcal Q_{\sigma}$-unique?
\end{problem}
Further, we want to address this issue thoroughly:
\begin{problem}
Can we provide a full characterization for a poset to be $\mathcal Q_{\sigma}$-unique?
\end{problem}

\end{document}